\documentclass[11pt]{article}

\usepackage{amsmath,amsfonts,amscd, amssymb,theorem}

\newcommand{\doubleheaddownarrow}{\big\downarrow\kern-3.325mm\downarrow}
\newcommand\Oh{\mathcal O}
\newcommand\Z{\mathbb Z}

\newcommand\om{\omega}

\newcommand{\iso}{\cong}

\newcommand\val{\operatorname{val}}

\newcommand\coker{\operatorname{coker}}
\newcommand\codim{\operatorname{codim}}
\newcommand\hatfunction{\operatorname{hat}}
\newcommand\grade{\operatorname{grade}}

\newcommand\Ext{\operatorname{Ext}}

\newcommand\Hom{\operatorname{Hom}}
\newcommand\res{\operatorname{res}}

\newcommand{\skewentry}{ \kern -1cm -\mathrm{sym} \kern -1cm}
\newcommand\Spec{\operatorname{Spec}}

 \newtheorem{theorem}{Theorem}[section]
 \newtheorem{lemma}[theorem]{Lemma}
 \newtheorem{prop}[theorem]{Proposition}
 \newtheorem{cor}[theorem]{Corollary}
 {
 \theorembodyfont{\rmfamily}
 \newtheorem{defn}[theorem]{Definition}

 \newtheorem{rem}[theorem]{Remark}

 }
 \newenvironment{pf}{\paragraph{Proof}}{\par\medskip}
 
 \newcommand{\qed}{\ifhmode\unskip\nobreak\fi\quad\ensuremath\square}
\newcommand{\QED}{\ifhmode\unskip\nobreak\fi\quad\ensuremath{\mathrm{QED}}}

\numberwithin{equation}{section}

\title{ Type II Unprojection }
 
\author{Stavros Argyrios Papadakis}
\date {January 2005}

\begin{document}
\maketitle

\begin {abstract}
Answering a question of M. Reid, we define and prove the
Gorensteiness of the type II unprojection. 
\end{abstract}

\section  {Introduction}

Unprojection is a philosophy, due to M. Reid, which is based on adjunction, 
and aims  to construct and analyse Gorenstein rings in terms of simpler 
ones. Geometrically, it is an inverse of projection and can also 
be considered as a modern version of Castelnuovo blow--down. 

The case that has been well understood by now is the unprojection of
type Kustin--Miller (or type I), which was studied
in \cite{PR} and \cite{P}, and is originally due, in a different 
formulation  using complexes, to  A. Kustin and M. Miller \cite{KM}. 
Examples by Reid,  A. Corti, S. Alt{\i}nok, G. Brown
and others (cf. \cite{Al}, \cite{CPR}, \cite{Ki}, \cite{R}) suggest that 
there  are also other types  of  unprojections that appear naturally 
and are useful in algebraic geometry.

The purpose of the present work is to provide the foundations
for the unprojection of type II. The initial data for the unprojection
is due to Reid, who also gave some indications about how it should
look like (cf. \cite {Ki} Section~9). The main difficulties 
compared to the unprojection of type Kustin--Miller are the existence 
of nonlinear 
relations and the increase of codimension by more that one. 
We give for the first time a precise definition 
(Definition~\ref{dfn!maindfn1}) of the generic type II 
unprojection and in Theorem~\ref{thr!gorensteiness}  we
prove that it is Gorenstein.

As a demonstration of our methods, we define using base change, and 
prove using generic perfection, the  Gorensteiness of the complete  
intersection type II unprojection (Section~\ref{sec!citypeII}). In  
Section~\ref{sec!calcandappl} we discuss calculational
aspects and 
applications to algebraic geometry. Finally, in Section~\ref{sec!finalcom}
we state some remarks and open questions.

ACKNOWLEDGEMENTS: I wish to thank Miles Reid for suggesting the problem,
encouragement, and many useful conversations, and Janko Boehm, Gavin
Brown,  Christoph Lossen and Frank--Olaf Schreyer for important discussions.
This work was financially supported by 
the Deutschen Forschungsgemeinschaft Schr 307/4-2.

\section {Generic type II unprojection}

\subsection {Definion of the unprojection ring}  \label{subs!dfnofunp}

Fix positive integers $k,n$ with $k \geq 1$ and $n \geq 2$. Our ambient ring
will be a polynomial ring $\Oh_{amb}$ over the integers
\[
           \Oh_{amb} = \Z [a_{ij}, z, w_p]
\]
where $1 \leq i \leq k+1$, $1 \leq j \leq n$ and $w_p$ is a finite number of 
indeterminants $w_1,w_2,  \dots $ that we use only for extra flexibility in 
the applications.  We assign $a_{ij}$ to have weight equal to $i$ for 
$1 \leq i \leq k+1$, $z$ to have weight $k+2$. We also give each $w_p$  
a positive number as weight.

Following Reid (cf. \cite{Ki} Section~9.5), we  define $I_D$  to be the ideal 
of $\Oh_{amb}$ generated by the $2 \times 2$ minors of the $2 \times n(k+1)$ 
matrix $M$, where 
\begin{equation}     \label{eq!dfnofM}
    M =  \begin {pmatrix}
             a_{21} & \dots & a_{2n}  &  \dots &
            a_{k+1,1} & \dots & a_{k+1,n}  &  z a_{11} & \dots & za_{1n}  \\
          a_{11} & \dots & a_{1n}  &  \dots &
            a_{k1} & \dots & a_{kn}  &  a_{k+1,1} & \dots & a_{k+1,n}\\      
     \end {pmatrix} 
\end{equation}
We notice that $I_D$ is a homogeneous ideal of $\Oh_{amb}$.  We set  
\[
        D = V(I_D) \subset \Spec \Oh_{amb},  \quad \Oh_D = \Oh_{amb} / I_D
\] 
and
\[
    \Oh_{\widetilde{D}} = \Z [x_{1},\dots ,x_{n}, t, w_p], \quad 
              \widetilde{D} = \Spec (\Oh_{\widetilde{D}}).
\]

\begin{rem} \label {rem!no1aboutD}
The ring $\Oh_D$  is not normal. Geometrically, the normalisation of $D$ is the map
\[
    q: \widetilde{D}  \to D 
\]
with
\[
       a_{ij} = x_jt^{i-1}, \quad z = t^{k+1}
\] 
for  $1 \leq i \leq k+1$, $1 \leq j \leq n$.  Since $q$ is an isomorphism in
codimension one, we have
\begin {equation}   \label{eq!normaliz}
   \om_D = q_* \om_{\widetilde{D}} = q_* \Oh_{\widetilde{D}}.
\end {equation}
\end {rem}
We also see that the ideal $I_D$ is prime of codimension $nk$ in $\Oh_{amb}$.

\paragraph {}
Let $I_X$ be a homogeneous prime ideal of $\Oh_{amb}$  such that $I_X \subset I_D$,
$I_X$ has codimemsion $nk-1$ in $\Oh_{amb}$,  and the integral domain 
$\Oh_X = \Oh_{amb}/I_X$ is normal and Gorenstein. We denote by $K(X)$ to be the
field of fractions of $\Oh_X$, and we set
\[
     I_D^{-1} = \{ f \in K(X) \colon f I_D \subset \Oh_X \}.
\]
We will consider $\Oh_D$ as a quotient ring of $\Oh_X$. Using 
(\ref{eq!normaliz}), $\om_D$ is a Cohen--Macaulay $\Oh_X$-module
which   needs $k+1$ generators 
$e_0, \dots e_k,$ as $\Oh_D$-module.
It is also clear  (cf.  \cite{Ki}~Section~9.5)
that we can choose them so that the module of relations between 
them is generated by relations of the form
\[
     M^t     
       \begin {pmatrix}
                e_i  \\
                e_{i+1}
     \end {pmatrix} = 0.
\]

\begin {defn}   \label {dfn!maindfn1}
The generic type II unprojection is  the  $\Oh_X$-subalgebra of $K(X)$
\[
         \Oh_X [I_D^{-1}] \subset K(X)
\] 
generated by $I_D^{-1}$.
\end {defn}
 The relation with adjunction theory is as follows:
 
Since, by assumption, $I_D$ has codimension one in the Gorenstein
ring $\Oh_X$, we have as in \cite{PR}~p.~563
a short exact sequence
\[
    0\to\om_X\to\Hom_{\Oh_X}(I_D,\om_X)\xrightarrow{\,\res_D\,}\om_D\to0
\]
where $\res_D$ is the  Poincar\'e residue map. 
Using the Gorensteiness of $\Oh_X$ we get the short exact sequence 
\begin{equation} \label{eq!adjunction}
   0\to\Oh_X\to\Hom_{\Oh_X}(I_D,\Oh_X)\xrightarrow{\,\res_D\,}\om_D\to0
\end{equation}
$I_D$ contains a regular element of $\Oh_X$, therefore 
$\Hom_{\Oh_X}(I_D,\Oh_X)$ is canonically isomorphic to $I_D^{-1}$,
and  we will  consider it as a submodule of $K(X)$.

\paragraph {}
After giving a presentation of $\Oh_X [I_D^{-1}]$  
(Proposition~\ref{prop!presentation}), we prove  
in Theorem~\ref{thr!gorensteiness} that it is Gorenstein.

\begin{rem}  \label{rem!notprincipal}
Assume the ideal $I_D \subset \Oh_X$ is principal, say $I_D = (a)$.
Then it clear that the ring   $\Oh_X [I_D^{-1}]$ is a
localization  of $\Oh_X$ hence Gorenstein. In the following, we
will assume that $I_D$ is not a principal ideal of $\Oh_X$.
\end{rem}

\subsection {Presentation of the unprojection ring}  \label{subs!presentat}

By the assumptions on $I_D$ and $\Oh_X$, we have a natural valuation map
\[
     \val_D \colon K(X)^* \to \Z 
\]
It is convenient to extend it by assigning $\val_D 0 = + \infty$.

\begin  {rem} \label {rem!no2aboutDE} 
For an element $s \in I_D^{-1}$ we have $\val_D (s) \geq -1$. 
We claim that every  $s \in \Oh_X [I_D^{-1}]$
with $\val_D (s) \geq 0$ is in $\Oh_X$. Indeed, by the normality
of $\Oh_X$ it is enough to prove that for every codimension one
prime ideal $I_E \subset \Oh_X$ we have $\val_E (s) \geq 0$. 
Let $I_E \not= I_D$, then there exists a point $p \in E$ and 
$f \in I_D$ with $f(p) \not= 0$. 
For sufficiently large $k$ we have $sf^k \in \Oh_X$, as
a consequence  $s$ is defined at the point $p$.  
\end{rem}

We fix $s_0, s_1, \dots ,s_k \in \Hom_{\Oh_X}(I_D,\Oh_X) \subset K(X)$ such 
that    $\res_D (s_i) = e_i$ for $0 \leq i \leq k$. By the adjunction 
sequence (\ref{eq!adjunction}), there exist elements  
$c_{i,j,p}, d_{j,p} \in \Oh_X$ 
(for $1 \leq i \leq k$, $1 \leq j \leq n$ and $0 \leq p \leq k-1$), 
such that 
\begin {eqnarray*}
     a_{i+1,j}s_p + a_{ij}s_{p+1} - c_{i,j,p}  & = &  0 \in K(X), \\
     z a_{1j}s_p + a_{k+1,j}s_{p+1}-d_{j,p}      & = &  0 \in K(X).
\end{eqnarray*}

Moreover, by the same sequence it follows that any (inhomogeneous) 
linear relation 
between $s_0, \dots ,s_k$ is generated by these relations.

Let $\Oh_X[T_0, \dots T_k]$ be the polynomial ring over $\Oh_X$ in
$k+1$ indeterminants, and denote by 
\[
     \phi \colon \Oh_X[T_0, \dots T_k] \to \Oh_X [I_D^{-1}]
\]
the $\Oh_X$-algebra homomorphism with  $\phi(T_i) = s_i$.
We set 
\begin {eqnarray*}
   f^a_{i,j,p}  & = &    a_{i+1,j}T_p + a_{ij}T_{p+1} - c_{i,j,p}, \\
   f^b_{j,p}    & = &     z a_{1j}T_p + a_{k+1,j}T_{p+1}-d_{j,p}.   
\end{eqnarray*}
These elements are in the kernel of $\phi$ and
moreover generate the submodule of linear relations. Our aim now
is to describe the quadratic relations and prove that together with 
the linear they generate the kernel of $\phi$. We need the following
lemma.

\begin {lemma}   \label{lem!quadraticrelations}
(Quadratic relations) The following are true.

a)  Let $ 0 \leq j,m \leq k, \quad 1 \leq i,l \leq k+1$ and     
     $1 \leq t \leq n$. Assume $l+m = i+j$. Then 
\[
        \val_D(a_{it}s_j - (-1)^{j-m}a_{lt}s_m) \geq 0. 
\]

b) Let $0 \leq i,j,l,m  \leq k$. Assume  $ l+m = i +j$. Then
\[
        \val_D(s_is_j - s_ls_m) \geq -1. 
\]

c)   Let $ 0 \leq i,j,l,m \leq k$.      
     Assume $l+m = i+j - (k+1)$. Then 
\[
        \val_D(s_is_j - (-1)^{k+1}zs_ls_m) \geq -1. 
\]   
\end{lemma}

\begin {pf}  
For a) we notice that assuming $i < l$, we have
\[
    a_{it}s_j - (-1)^{j-m}a_{lt}s_m 
    = v_i - v_{i+1} + v_{i+2} - \dots + (-1)^{l-i-1} v_{l-1},
\] 
where $ v_q = a_{qt}s_{j+i-q}+a_{q+1,t}s_{j+i-q-1}$.
Taking into account the linear relations the result follows.

For b) it is enough to show that 
 \[
         \val_D (s_is_j-s_{i-1}s_{j+1}) \geq -1. 
\]
We have
\[
    a_{11}(s_is_j-s_{i-1}s_{j+1}) = s_j(a_{11}s_i+a_{21}s_{i-1})
        -s_{i-1}  (a_{11}s_{j+1} +a_{21}s_j).
\]
Taking into account the linear relations the result follows.

For c), using b) it is enough to show that for $i \geq 1$ 
we have 
\[
         \val_D (s_ks_i-(-1)^{k+1}zs_{i-1}s_0) \geq -1. 
\]
We have 
\begin {eqnarray*}
 & & a_ {11}( s_ks_i-(-1)^{k+1}zs_{i-1}s_0 )  = \\
 & &    s_i (a_{11}s_k - (-1)^ka_{k+1,1} s_0)  -
    (-1)^{k+1} s_0 (a_{11}zs_{i-1} + a_{k+1,1}s_i), 
\end{eqnarray*}
so the result follows from a).
 \QED \medskip
\end{pf}

\paragraph{}
Using Proposition~\ref{lem!quadraticrelations}, 
for $i,j$ with $i+j \leq k$ there exists
polynomial $g^a_{ij}(T_0, \dots ,T_k) \in \ker(\phi)$,
of the form 
\[
    g^a_{ij} = T_iT_j - T_0T_{i+j} + \text {linear terms}.
\]
In addition, for $i,j$ with $i+j \geq k+1$,  there exist
polynomial $g^b_{ij}(T_0, \dots ,T_k) \in \ker(\phi)$,
of the form 
\[
   g^b_{ij} = T_iT_j - (-1)^{k+1}zT_0T_{i+j-k-1} + \text {linear terms}.
\]
In both cases by linear terms we mean terms of total degree in $T_0, 
\dots T_n$ at most one.

\begin {prop}   \label{prop!presentation}
We have
\[
      \ker (\phi) = ( f^a_{i,j,p}, f^b_{j,p},g^a_{ij}, g^b_{ij}),
\]
with indices as above.
\end{prop}

\begin{pf}

Taking into account  the relations $g^a,g^b$, it is enough to restrict to  polynomials 
$h(T_0,\dots ,T_k)$ which have  total degree with respect to the variables 
$(T_1, \dots ,T_k)$ at most one.

We use induction on the total degree (with respect to all variables) $p$ of the
polynomial  $h$.

We already noticed  that for $p=1$ it is true.
   
Assume it is true for all polynomials of total degree at most $p-1$ in the kernel of $\phi$  
which  have  total degree with respect to the variables $(T_1, \dots ,T_k)$
at most one.  Let $h(T_0,T_1, \dots ,T_k)$ be a polynomial 
of total degree exactly $p$ in the kernel of $\phi$, which has $(T_1, \dots ,
T_k)$-degree at most one. Assume
\[
    h(T_0,\dots ,t_k) = e_0 T_0^p+T_0^{p-1}u_1(T_1,\dots,T_k)+ L(T_0, \dots ,T_k)
\]
where $e_0 \in \Oh_X$, $u_1$ is (homogenously) linear in $(T_1,\dots ,T_k)$, 
and $L(T_0, \dots ,T_k)$ has total degree at most $p-1$. Write 
\[
      h(T_0,\dots ,T_k)= T_0^{p-1}(e_0 T_0+ u_1)+  L.
\]
We have 
\[
    \val_D(e_0s_0+ u_1(s_1,\dots ,s_k)) \geq 0,
\]
otherwise $\val_Dh(s_0, \dots ,s_k) \leq -k$, which contradicts $h(s_0, \dots ,s_k) =0$. 
Therefore, by Remark~\ref{rem!no2aboutDE}  
\[
     e_0s_0+ u_1(s_1,\dots ,s_k)  \in \Oh_X.
\]
By the case $p=1$, there exists $c \in \Oh_X$ such that  
\[ 
     e_0T_0+ u_1(T_1,\dots ,T_k) -c  \in (f^a, f^b, g^a, g^b).
\]
We finish the proof by using the inductive hypothesis on the polynomial 
\[
    c T_0^{p-1}  + L(T_0, \dots ,T_k).
\]
\QED \medskip
\end{pf}

\subsection {The theorem}

\begin{rem}  \label{remark!no9}
As in \cite{PR}~Lemma~1.1,   we can assume  that
$s_0 \in \Hom_{\Oh_X}(I_D,\Oh_X)$ is injective, write
\[
  s_0 \colon I_D \to I_N \subset \Oh_X
\]
for the image of $s_0$. 
\end{rem}

We will use the following rather general result

\begin{prop} \label{prop!codimgenresult}
Let $A$ be a  Cohen--Macaulay ring, $I$ an
ideal of codimension one, $\phi \colon I \to A$ an injective
$A$-homomorphism with image the ideal $J \subset A$. Assume
that if for $a \in A$ we have $I \subset (a)$, then $a$ is 
a unit of $A$. Then $J$ has codimension one in $A$.
\end{prop}

\begin {pf}
By the assumption, there exists an $A$-regular element $w$ in $I$.
$\phi(w)$ is not a unit, otherwise for all $b \in I$ we have 
\[
  b = \phi(w) (\phi(w))^{-1}b = (\phi(w))^{-1}\phi(b))w,
\]  
contradicting
our assumptions. Since $\phi$ is injective, $\phi(w)$ is a regular
element of $J$, so $J$ has codimension at least one.

Assume now that $J$ has codimension at least two.  Since $A$
is Cohen--Macaulay, there exists an $A$-regular sequence of length
two contained in $J$. That immediately implies that the natural map
from $A$ to $\Hom_A(J,A)$ is an isomorphism. As a result,
$\phi^{-1} \colon J \to I$ is given by multiplication by an 
element $a \in A$. If $a$ is a unit, we get $J=I$, a contradiction
since $I$ has codimension one. If $a$ is not a unit, then 
$I \subset (a)$, contradicting our assumptions. So $J$ can only
have codimension one.
 \QED \medskip
\end{pf}

\begin {cor}  \label{cor!codimisone}
 We have $\codim_{\Oh_X}I_N = 1$.
\end {cor}

\begin {pf}   Suppose that for an element
$a \in \Oh_X$ we have
\[
           I_D \subset (a)  \subset \Oh_X.
\]
Recall that we assume (Remark~\ref{rem!notprincipal})
that the ideal $I_D$ of $\Oh_X$ is not principal, therefore
$a \notin I_D$. Let $w_1 \in I_D$ be a nonzero element. 
We have that $a$ 
divides $w_1$, and since $I_D$ is prime, there exists 
$w_2 \in I_D$ with $w_1 = aw_2$. Continuing this way,
we get an increasing sequence of ideals of $\Oh_X$
\[
    (w_1) \subset (w_2) \subset (w_3) \subset  \dots
\]
such that $w_i = a w_{i+1}$. Since $\Oh_X$ is Noetherian,
this sequence becomes stationary, so there exists an index $i$ and
$b \in \Oh_X$ with $w_{i+1} = b w_i$. Therefore 
$(ab-1)w_i =0$, and since $\Oh_X$ is an integral domain we get 
$a$ a unit.  Then, the Corollary follows from  
Proposition~\ref{prop!codimgenresult}.
 \QED \medskip
\end{pf}

\paragraph{}  For simplicity, we set 
\[
    R_1 =  \Oh_X[T_0,\dots ,T_k]/(f^a, f^b,g^a, g^b). 
\]
By Proposition~\ref{prop!presentation}, $R_1$ 
is isomorphic under $\phi$ to $ \Oh_X[I_D^{-1}] \subset K(X)$,
so it is an integral domain and therefore $T_0$ is an $R_1$-regular
element.  Let $(T_0)$ be the ideal of $R_1$ generated by $T_0$.

\begin {prop}  \label{prop!threeisom}
We have the following isomorphisms of $\Oh_X$-modules
\[
   R_1/(T_0)  \iso  \Hom_{\Oh_X}(I_D,\Oh_X)/ (s_0)  \iso 
         \Hom_{\Oh_X}(I_N,\Oh_X)/(i_N) \iso  \om_N,
\]
where $i_N \in \Hom_{\Oh_X}(I_N,\Oh_X)$ is the natural inclusion 
$ i_N \colon I_N \to \Oh_X$. Notice that in the second and third quotient 
we divide by the $\Oh_X$-submodules generated by the respective 
elements, while in the first with the ideal (that is, the $R_1$-submodule) 
generated by $T_0$.
\end{prop}    
    
\begin{pf} 
  By Corollary~\ref{cor!codimisone} we have that $N$ has codimension one in $X$.
 Therefore, the adjunction exact sequence 
\[
   0 \to \Oh_X \to  \Hom_{\Oh_X}(I_N,\Oh_X) \to \om_N \to 0
\]
gives the isomorphism 
\[
     \Hom_{\Oh_X}(I_N,\Oh_X)/(i_N) \iso  \om_N.
\]

The isomorphism $s_0 \colon I_D \to I_N$ induces an isomorphism
\[
   s_0^* : \Hom_{\Oh_X}(I_N,\Oh_X) \to \Hom_{\Oh_X}(I_D,\Oh_X)
\]
with  $(s_0^*)(i_N) = s_0$. The isomorphism 
\[
     \Hom_{\Oh_X}(I_D,\Oh_X)/ (s_0)  \iso \Hom_{\Oh_X}(I_N,\Oh_X)/(i_N)
\]
follows.

Now we will prove that
\[
     R_1/(T_0)  \iso  \Hom_{\Oh_X}(I_D,\Oh_X)/ (s_0).
\]
Using the isomorphism  $ R_1 \iso  \Oh_X[I_D^{-1}]$, is it enough to prove that
\[
   \Oh_X[I_D^{-1}]/ (s_0)  \iso \Hom_{\Oh_X}(I_D,\Oh_X)/ (s_0)
\]
as $\Oh_X$-modules. Again, the first quotient is with respect to the ideal generated by $s_0$, 
while the second is with respect to the submodule.

Denote by $\psi$ the composition 
\[
     \Hom_{\Oh_X}(I_D,\Oh_X) \to  \Oh_X[I_D^{-1}] \to  \Oh_X[I_D^{-1}]/(s_0),
\]
(the first map is the natural inclusion).

Clearly $s_0$ is in the kernel  of $\psi$.  Using the relations $g^a,g^b$
we get that the map $\psi$ is surjective.

\paragraph {CLAIM} $\ker (\psi) = (s_0)$.

Proof of Claim: Assume $u(s_1, \dots ,s_k) \in \ker(\psi)$, where  
$u(T_1, \dots ,T_k) \in \Oh_X[T_1, \dots ,T_k]$ has total degree at
most one. There exists
$h(T_0, \dots ,T_K) \in \Oh_X[T_0, \dots ,T_k]$ with
\[
    u(s_1, \dots ,s_k) = s_0 h(s_0, \dots ,s_k).   
\]
(Equality inside $\Oh_X[I_D^{-1}]$.) 
 As a consequence, $\val_D(h(s_0, \dots ,s_k)) \geq 0$. By 
Remark~\ref{rem!no2aboutDE}, $h(s_0,\dots ,s_k) \in \Oh_X$.
Therefore, $u(s_1, \dots ,s_k)$ is inside the $\Oh_X$-submodule 
of   $\Hom_{\Oh_X}(I_D,\Oh_X)$ generated by $s_0$.
\QED \medskip
\end{pf}

\begin{rem}  \label{remark!no10}
To prove that $R_1$ Gorenstein it is enough to prove that $R_1/(T_0)$ is 
Gorenstein.  This is because $I_D$ and $\Oh_X$ are graded so
$\Oh_X[I_D^{-1}]$ has a natural  grading and  with that,  $T_0$ is  
a homogeneous element of positive degree 
(cf. \cite{BH} Exerc.~3.6.20(c), p.~144).  
\end{rem}

\begin{rem}  \label{remark!no11}
To prove that $R_1/(T_0)$ is Gorenstein, it is enough, due to 
Proposition~\ref{prop!threeisom}, to prove 
\[
     \grade_{\Oh_X}(m, \om_N) = \dim \Oh_X -1, 
\]
where $m = (a_{ij},z,w_p)$ is the unique maximal homogeneous ideal 
of $\Oh_X$, and
\[
    \Ext^1_{\Oh_X}(\om_N, \Oh_X) \iso \om_N,
\]
the isomorphism as $\Oh_X$-modules.
\end{rem}

\begin {prop}  \label{prop!depthOK}
We have 
\[
     \grade_{\Oh_X}(m, \om_N)  = \dim \Oh_X -1.
\]
\end {prop}

\begin{pf}   By Remark~\ref{rem!no1aboutD}, 
  \[
      \grade_{\Oh_X}(m, \om_D)  = \dim \Oh_X -1.
 \]
 Using the two adjunction sequences for $\om_D$ and 
 $\om_N$ the result follows using the arguments
of  step 2 of  the proof of  \cite{PR}  Theorem~1.5. 
\QED \medskip
\end{pf}

\begin{rem}  \label{remark!no12}
  It is clear that the map 
\[
    h_3  \colon \Oh_D \to \om_D 
\]
with
\[
    h_3(a) = a \res(s_0)
\] 
is injective. Write $T$ for the cokernel of $h_3$.  From Remark~\ref{rem!no1aboutD}
it follows that as an $\Oh_X$-module $T$ is isomorphic  to a direct sum of
$k$ copies of $\Oh_X / (a_{11}, \dots , a_{k+1,n})$.

Therefore, the support of $T$ in $\Oh_X$ has codimension in $\Spec(\Oh_X)$ 
equal to $(k+1)n- (nk-1) = n+1 \geq 3$. Using \cite{BH} Corollary~3.5.11, 
we get 
\[
     \Ext^1_{\Oh_X}(T,\Oh_X) =  \Ext^2_{\Oh_X}(T,\Oh_X) = 0.
\]
\end{rem}

\paragraph {} We are now ready to prove our main theorem.

\begin {theorem}   \label{thr!gorensteiness}
   The ring  $\Oh_X[I_D^{-1}]$ is Gorenstein 
\end{theorem}

\begin{pf}  
We follow the pattern of the proof of \cite{PR} Theorem~1.5. 

Using Remarks~\ref{remark!no10},~\ref{remark!no11}, and 
Proposition~\ref{prop!depthOK}, it is enough to prove that 
\[
    \Ext^1_{\Oh_X}(\om_N, \Oh_X) \iso \om_N
\]
as $\Oh_X$-modules. 
 
Start with the two columns consisting of canonical maps 
\[
 \renewcommand{\arraystretch}{1.5}
 \begin{matrix}
 &&I&&\Oh_X&&&& \\
 && \bigcap && \bigcap \\
 &&\Oh_X&&\Hom_{\Oh_X}(I_D,\om_X)&&& \\
 && \doubleheaddownarrow && \doubleheaddownarrow \\
 &&\Oh_D&&\om_D \\
 \end{matrix}
 \]
We define three maps 
\begin {eqnarray*}
   h_1 &\colon&  I_D \to \Oh_X,  \quad  h_1(a) = s_0(a) \in \Oh_X, \\
   h_2 &\colon&  \Oh_X \to  \Hom_{\Oh_X}(I_D,\om_X), \quad (h_2(a))(b)=s_0(ab) \in \Oh_X, \\
   h_3 &\colon&  \Oh_D \to \om_D, \quad h_3(a) = a \res(s_0) \in \om_D.
\end{eqnarray*}

It is clear that all three maps $h_1, h_2, h_3$ are injective (compare Remark~\ref{remark!no12}).
In addition, 
\[
\coker h_1 = \Oh_N, \quad \coker h_2 \iso \om_N, \quad  \coker h_3 = T,
\]
where $T$ was calculated in Remark~\ref{remark!no12}.   We have a commutative diagram

\[
 \renewcommand{\arraystretch}{1.5}
 \begin{matrix}
 0&\to&I&\xrightarrow{\ h_1\ }&\Oh_X&\to&\Oh_N&\to&0 \\
 && \bigcap && \bigcap \\
 0&\to&\Oh_X&\xrightarrow{\ h_2\ }&\Hom_{\Oh_X}(I_D,\Oh_X)&\to&\om_N &\to&0 \\
 && \doubleheaddownarrow && \doubleheaddownarrow \\
 0&\to&\Oh_D&\xrightarrow{\ h_3\ }&\om_D &\to& T &\to& 0\\
 \end{matrix}
 \]
which induces a map $d \colon \Oh_N \to \om_N$. We claim that this map is injective.
Indeed, let $x \in \Oh_X$ such that $\res_N \circ s^* (m_x) = 0 \in \om_N$, where 
$m_x \in \Hom_{\Oh_X}(I_D,\Oh_X)$ is the multiplication by $x$. 
Then there exists $w \in \Oh_X$ such that 
\[
     m_x = h_2(w) 
\]
hence,
\[
    xa = ws_0(a)   \in \Oh_X
\] 
for all $a \in I_D$, which implies that  
\[
   x-ws_0 =0 \in  \Hom_{\Oh_X}(I_D,\om_X).
\]
Using the valuation $\val_D$ we get $w \in I_D$. Therefore
$x = h_1(w)$.

\paragraph{} Using the Snake Lemma, we get an exact sequence
\[
    0 \to \Oh_N \to \om_N \to T \to 0
\] 
Taking the long exact associated to $\Hom_{\Oh_X}( -,\om_X)$,
and using the vanishings from Remark~\ref{remark!no12}, we get
\[
   \Ext^1_{\Oh_X}(\om_N, \Oh_X) \iso  \Ext^1_{\Oh_X}(\Oh_N, \Oh_X) \iso  \om_N,
\]
which finishes the proof of the theorem.
\QED \medskip
\end{pf}

The next proposition is a corollary of Proposition~\ref {prop!presentation}
and  Theorem~\ref{thr!gorensteiness}.

\begin {prop}   \label{prop!presentationno2}
The map $\phi$ induces an isomorphism
\[
    \Oh_X[I_D^{-1}]  \iso  \frac {\Oh_{amb}[T_0, \dots ,T_k]}{I_X + (f^a,f^b,g^a,g^b)}.
\]
Moreover,
if we consider $\Oh_X[I_D^{-1}]$ as an  $\Oh_{amb}[T_0, \dots ,T_k]$-module
via $\phi$, it is perfect of grade $(nk-1)+(k+1)=nk+k$.  
\end{prop}

\begin {pf} 
The isomorphism is immediate from  Proposition~\ref {prop!presentation}.
Using \cite{BV} p.~210 remarks after Proposition~16.19, the  perfectness
follows from Theorem~\ref{thr!gorensteiness}. We have
\[
   \grade \Oh_X[I_D^{-1}] = \dim \Oh_{amb}[T_0, \dots ,T_k] - \dim  \Oh_X[I_D^{-1}].
\]
Hence, it is enough to prove that $\dim \Oh_X  =  \dim \Oh_X[I_D^{-1}]$. Set 
\[
       J = (f^a,f^b,g^a,g^b) \cap \Oh_X[T_0].
\]
 Using the relations $g^a,g^b$ we see that  $\Oh_X[I_D^{-1}]$ is a finitely
generated   $\Oh_X[T_0]/ J$-module. As a result, the extension 
\[
      \Oh_X[T_0]/ J \subset   \Oh_X[I_D^{-1}]
\]
is integral, so it is enough to prove $\dim \Oh_X  = \dim \Oh_X[T_0]/ J$.
This follows from the facts that $T_0$ is a regular element of  
$\Oh_X[T_0]/ J$ and the quotient by the ideal generated by $T_0$ is
isomorphic to the ring $\Oh_X / I_N$ (cf. proof of Theorem~\ref{thr!gorensteiness}).
\QED \medskip
\end{pf}

\section {Complete intersection type II unprojection}   \label{sec!citypeII}

\subsection {Generic complete intersection type II}  \label{subs!genci}

The generic complete intersection type II is the specialization of the generic type
II to the case where $I_X$ is a complete intersection. More precisely,
the initial data is the triple 
\[
    I^g_X \subset I^g_D \subset \Oh^g_{amb}
\]
(g for generic) defined as follows. 

Fix as in Subsection~\ref{subs!dfnofunp} positive integers $k,n$ with $k \geq 1$ 
and $n \geq 2$.
$I^g_D $ is the ideal generated by the $2 \times 2$ minors $u_1, \dots ,u_l$ 
of the matrix $M$ defined in (\ref{eq!dfnofM}). We define
\[
      I^g_X = (f^1, \dots ,f^{nk-1}), 
\]
where for $1 \leq  p \leq  nk-1$ we have 
\[
    f^p = \sum_{j=1}^l w_{pj}u_{j}.
\]
Finally $\Oh^g_{amb}$, is the polynomial ring over the integers in all the
indeterminants appearing above,
\[
    \Oh^g_{amb} = \Z [a_{ij}, z, w_{pj}].
\]
We put weights for $a_{ij}$ and $z$ as in  Subsection~\ref{subs!dfnofunp},
while for $w_{pj}$ we make the unique minimal choice of positive weights
such that all $f^p$ become homogeneous.

Our configuration satisfies the conditions of Subsection~\ref{subs!dfnofunp},
so setting, for simplicity,
\[
    R_1^g  = \Oh_X^g [(I_D^g)^{-1}],
\]
we have that $R_1^g$ is a Gorenstein ring (Theorem~\ref{thr!gorensteiness}).
In addition it has the presentation described in 
Proposition~\ref{prop!presentationno2}  and by the same Proposition  
is a perfect  $\Oh_{amb}^g[T_0, \dots ,T_k]$-module of grade $nk+k$.

\subsection {Complete intersection type II}

\paragraph {} 
Assume we have a triple 
\begin {equation} \label{eq!cidata}
    I^s_X \subset I^s_D \subset \Oh^s_{amb}
\end{equation}
(s for specific) as follows. 

$I^s_D $ is the ideal generated by the $2 \times 2$ minors 
$\widehat{u_1}, \dots ,\widehat{u_l}$  of a matrix $\widehat{M}$,
of the form defined  in (\ref{eq!dfnofM}), but with
elements $\widehat{a_{ij}}$ instead of $a_{ij}$, and $\widehat{z}$
instead of $z$. We assume that $I^s_D$ has codimension in
$\Oh^s_{amb}$ equal to $nk$,  which is the maximal possible.

$\Oh^s_{amb}$ is an equidimensional Gorenstein ring, not necessarily 
local or graded, containing all elements with hats.

$I^s_X$ has codimension $nk-1$ in $\Oh^s_{amb}$ and is generated
by an $\Oh^s_{amb}$-regular sequence $\widehat{f^1}, \dots ,
\widehat{f^{nk-1}}$, so there are $\widehat{w_{pj}} \in \Oh^s_{amb}$ 
such that   
\[
   \widehat{f^p} = \sum_{j=1}^l \widehat{w_{pj}}\widehat{u_{j}}.
\]
for  $1 \leq  p \leq  nk-1$.

\paragraph { }We consider the ring homomorphism 
\[
     \hatfunction \colon  \Oh^g_{amb} [T_0, \dots ,T_k] \to  
         \Oh^s_{amb} [T_0, \dots ,T_k]
\]
which sends elements without hats to the corresponding elements 
with  hats and $T_i$ to itself.

\begin {defn}   \label {dfn!cimaindfn}
The type II unprojection  for our initial data (\ref{eq!cidata})
is the $\Oh_X^s$-algebra
\[
        R_1^s = R_1^g \otimes  \Oh^s_{amb} [T_0, \dots ,T_k],
\]
the tensor product over the ring $\Oh^g_{amb} [T_0, \dots ,T_k]$.
\end {defn}

\paragraph { }We have the following theorem.

\begin{theorem}  
    The ring $R_1^s$ is Gorenstein.
\end{theorem}

\begin {pf}
   By  \cite{BV}~Theorem~3.5,  it is enough to prove that 
the ring $R_1^s$ has grade at least equal to $nk+k$ as  
$\Oh^s_{amb} [T_0, \dots ,T_k]$-module.  For that  we argue as in  
Proposition~\ref{prop!presentationno2}, taking into 
account that $I_D^s$ contains an $\Oh_X^s$-regular element.
\QED \medskip
\end{pf}

\begin{rem}  \label{remark!no1001}
The reason we defined the unprojection using base change   
is that we expect that the ring $\Oh_X^s [(I_D^s)^{-1}]$ 
is not always Gorenstein. An example in the
unprojection of type Kustin--Miller is  the following. Let $R$ be
the homogeneous coordinate ring  of the plane cuspidal cubic 
$x_2^2x_0-x_1^3$ and $I=(x_1,x_2)$ the reduced ideal of
the cusp. Then the ring $R[I^{-1}]$ is isomorphic
to the  homogeneous coordinate ring of the twisted cubic, 
hence it is not Gorenstein. 
\end{rem}

\section {Calculations and Applications}  \label{sec!calcandappl}

\subsection {Computations}

\paragraph{} The simplest case, due to Reid,  is the complete intersection 
type II unprojection for parameters $k=1$ and $n=2$, in which we pass from 
codimension
$1$ to a codimension $3$ Pfaffian. We refer to \cite{Ki}~Section~9.5 for 
the explicit calculations.

\paragraph {} In $\cite{P2}$, which is work under progress, we explicitly 
calculate  the linear relations for the type II  unprojection for parameters 
$k=1$ and any $n \geq 2$. The main idea is to identify the
projective resolution of $\om_D$ as a twisted sum of two Koszul complexes.

The calculation of the quadratic relations turns out to be more difficult.
We expect that, assuming $2$ is invertible, there are symmetric
formulas, but so far, we have been able to write them down only for
the complete intersection type II with $k=1$ and $n  \leq 3$.

\subsection {Applications}

One of the first appearances of the type II unprojection  in algebraic geometry 
is in the elliptic involutions between Fano threefold hypersurfaces in 
\cite{CPR}. According to  \cite{CPR}~Section~7.3, it is expected 
that the the various unprojection types  will eventually serve as the basis 
for a  classification for all (possibly singular) Fano threefolds.

Another application of the type II unprojection, is in the construction
of new K$3$ surfaces, Fano threefolds, and Calabi--Yau threefolds.
In \cite{Al}, Alt{\i}nok constructs a large number of codimension 
three K$3$ surfaces using type II unprojection, and conjectures the 
existence of more K$3$ surfaces and Fano threefolds, including a 
codimension four Fano with $H^0(-K_X) = 0$ (cf. \cite {Ki} Example~9.14).

In the Magma K3 database (cf. \cite {ABR}~Section~6)  due to Brown,
there is a large number of
additional K$3$ surfaces which are likely to exist as
type II unprojections, including cases in  codimensions 
$6,7$ and $8$. We believe that the ideas of the present work together
with explicit singularity calculations can  establish the existence of 
many of them.

Finally, there is work in progress by A. Buckley and  B. Szendr{\H{o}}i, 
which is expected to  lead to the construction of  new Calabi--Yau 
threefolds using type II unprojection (cf. \cite{BS}).

\section {Final remarks and questions}  \label{sec!finalcom}

\begin{rem} \label {remark!no333}

An  interesting 
problem is to find conditions which will guarantee 
that $R_1^s = \Oh_X^s [(I_D^s)^{-1}]$, compare Remark~\ref{remark!no1001}.
More generally, 
for given unprojection data, how can we intrinsically 
distinguish  the relations of 
$\Oh_X^s [(I_D^s)^{-1}]$ that are also
present in the generic case from those that
occur because $\Oh_X^s$ is not sufficiently generic?
\end{rem}

\begin {rem} \label{remark!no334}
How can we construct a projective resolution of $R_1$ as 
\\
$\Oh_{amb}[T_0, \dots ,T_k]$-module using projective
resolutions of $\Oh_D$ and  $\Oh_X$?  In the case 
of unprojection of type Kustin--Miller it was done
in \cite{KM}.
\end{rem}

\begin {rem} \label{remark!no335} 
What other types of type II unprojection format exist, in
addition to the complete intersection one?
\end {rem}

\begin{rem} \label{remark!no336} 
In \cite{CR}, Corti and Reid pose the problem of interpreting
the Gorenstein formats arising from unprojection as 
solutions to universal problems. What can be said about
the type II case?
\end{rem}

\bigskip
\noindent
Stavros Papadakis, \\
Mathematik und Informatik,  Geb. 27 \\
Universitaet des Saarlandes \\
D-66123 Saarbruecken, Germany \\
e-mail: spapad@maths.warwick.ac.uk


\begin {thebibliography} {xxx}      


\bibitem[Al]{Al}
Alt{\i}nok S.,
\textsl{Graded rings corresponding to polarised
K3 surfaces and $\mathbb Q$-Fano 3-folds}.
Univ. of Warwick Ph.D. thesis,
Sept. 1998, 93+ vii pp.


\bibitem[ABR]{ABR}
Alt\i nok, S.,  Brown, G.  and Reid, M., \textsl {
Fano 3-folds, $K3$ surfaces and graded ring}, in 
Topology and geometry: commemorating SISTAG, 
Contemp. Math., 314,
AMS 2002, 25--53,



\bibitem[BH]{BH} Bruns, W. and Herzog, J.,  \textsl{ 
Cohen-Macaulay rings}. 
Revised edition, 
Cambridge Studies in Advanced Mathematics 39, CUP 1998

\bibitem[BV]{BV}
Bruns, W. and Vetter, U.,  \textsl{ 
Determinantal rings}. 
Lecture Notes in Math. 1327, 
Springer 1988

\bibitem[BS]{BS} Buckley, A. and Szendr{\H{o}}i, B., \textsl{ Orbifold 
Riemann--Roch for threefolds with an application to Calabi--Yau geometry}, 
submitted, math.AG/0309033, 19~pp.



\bibitem[CR]{CR}
Corti A. and Reid M., \textsl {
Weighted Grassmannians}, in 
Algebraic geometry, A volume in memory of Paolo 
Francia, 
M. Beltrametti et al (eds.), de Gruyter 2002, 141--163


\bibitem[CPR]{CPR}
Corti A., Pukhlikov A. and Reid M.,
\textsl{Birationally
rigid Fano hypersurfaces}, in Explicit birational geometry  
of 3-folds, 
A. Corti and M. Reid (eds.), CUP 2000, 175--258

\bibitem[KM]{KM} Kustin, A. and Miller, M., \textsl { 
Constructing big Gorenstein ideals from small ones}. 
J. Algebra {\bf 85} (1983),  303--322

\bibitem[P]{P} Papadakis, S., \textsl {
Kustin--Miller unprojection with complexes},
J. Algebraic Geometry {\bf 13}  (2004), 249-268 

\bibitem[P2]{P2} Papadakis, S., \textsl { The relations of type II
unprojection}, work in progress

\bibitem[PR]{PR} Papadakis, S. and Reid, M., \textsl {
Kustin--Miller unprojection without complexes},
J. Algebraic Geometry {\bf 13}  (2004), 563-577 


\bibitem[R]{R}
Reid, M.,
\textsl{Examples of type IV unprojection},
math.AG/0108037, 16~pp.

 \bibitem[Ki]{Ki} Reid, M., \textsl {
Graded Rings and Birational Geometry}, 
in Proc. of algebraic symposium (Kinosaki, Oct 2000), 
K. Ohno (Ed.) 1--72, available from  
www.maths.warwick.ac.uk/$\sim$miles/3folds


\end{thebibliography}
\end{document}